\documentstyle[12pt]{article}

\begin{document}

\author{ T. Kadeishvili }
\title{On the Differentials of the Spectral Sequence of a Fibre Bundle}
\date{English translation of the paper from Bulletin of the Academy of Sciences of the Georgian SSR, 82 N 2, 1976}
\maketitle

The theorem of W. Shi, \cite{shih}, which in its turn generalizes
the theorem of Faddell-Hurewicz \cite{fadel}, states the
following: if the structure group $G$ of a fibre bundle $F\to E\to
B$ is $(n-1)$-connected, then in the spectral sequence of this
bundle the differentials $d^r=0$ for $r<n$ and for $n\leq r\leq
2n-2$ are determined by certain characteristic classes of the
associated principal bundle. The above mentioned Faddell-Hurewicz
theorem calculates only the first nontrivial differential $d^n$.

For simplicity, to avoid difficulties with signees, all homologies
bellow we consider over $Z_2$ although it is enough to assume them
free.

In this paper we generalize the Shis result for higher dimensions:
we construct some cochains from $C^*(B,H_*(G))$ and certain {\it
polylinear operations} which define all differentials $d^r$. In
dimensions $n\leq r\leq 2n-2$ these cochains are cocycles form the
classes from Shis theorem and the above mentioned polylinear
operations is just $\smile$ product.

\vspace{5mm} In papers \cite{ber1}, \cite{ber2} N.A. Berikashvili
has constructed cochain $h=h^2+h^3+...,\ h^r\in
C^r(B,Hom(H_*(F),H_*(F))$ (representative of the {\it
predifferential} of the fibration) which determines all
differentials of the spectral sequence. From this point of view
the Shis theorem states that for $(n-1)$-connected structure group
one has $h^r=0$ for $r<n$ and for $n\leq r\leq 2n-2$ the
components $h^r$ can be expressed in terms of some cocycles from
$C^r(B,H_{r-1}(G))$.

\vspace{5mm} Let $\xi=(X,p,B,G)$ be a principal $G$-bundle, $F$ be
a $G$ space and $\eta=(E,p,B,F)$ be the associated bundle with the
fiber $F$.

\vspace{5mm} Let $C^*=C^*(H_*(G),H_*(G))$ be the Hochschil cochain
complex of the ring $H_*(G)$ with coefficients in itself: cochains
in $C^i$ are homomorphisms $f^i:H_*(G)\otimes
...(i-times)...\otimes H_*(G)\to H_*(G)$ and the coboundary
operator is given by $\delta f^i(a_1,...,a_{i+1})=a_1\cdot
f^i(a_2,...,a_{i+1})+\sum_k f^i(a_1,...,a_{k}\cdot
a_{k+1},...,a_{i+1})+f^i(a_1,...,a_{i})\cdot a_{i+1}$.

There are $\smile$ and $\smile_1$ products in the Hochschild
complex $C^*$ defined as follows: for $f\in C^i$ and $g\in C^j$
let
$$
\begin{array}{c}
f\smile g(a_1,...,a_{i+j})=f(a_1,...,a_i)\cdot
g(a_{i+1},...,a_{i+j});\\
f\smile_1 g(a_1,...,a_{i+j-1})=\sum_k f(a_1,...,a_k,
g(a_{k+1},...,a_{k+j}),a_{k+j+1},...,a_{i+j-1}).
\end{array}
$$

The standard formulas $\delta(f\smile g)=\delta f\smile g+f\smile
\delta g,\ \ \delta(f\smile_1 g)=\delta f\smile_1 g+f\smile_1
\delta g +f\smile g+g\smile f$ are valid for so defined $\smile$
and $\smile_1$.

\vspace{5mm} Now let $\bar{C^*}=C^*(C^*(B,H_*(G)),C^*(B,H_*(G))$
be the Hochschild complex of the ring $C^*(B,H_*(G))$. There
exists a map $\mu:C^*\to \bar{C^*}$ which assignees to $f^i\in
C^i$ the following homomorphism: for $b^k\in C^*(B,H_*(G))$ let
$(\mu f)(b_1,...,b_i)$ be the $i$ fold $\smile$ product of
elements $b_1,...,b_i$ when the coefficients are multiplied by
$f^i$.

\vspace{5mm} Let us define a {\it Hochschild twisting cochain } as
an element $f=f^3+f^4+...,\ f^k\in C^k$ satisfying the condition
$\delta f=f\smile_1f$. The set of all such twisting cochains we
denote by $N$.

Let us define the subset $L\subset N\times C^*(B,H_*(G))$ as
$L=\{(f,h_0),\ f\in N,\ h_0\in C^*(B,H_*(G)),\ s.t. \ \delta
h_0=h_0\cdot h_0+\sum_k (\mu f^k)(h_0,...,h_0)\}$.

\vspace{5mm} {\bf Definition of the map $\alpha:L\to
D^*(B,H_*(G))$.} First recall \cite{ber1}, \cite{ber2} that
$D^*(B,H_*(G))$ is defined as the factorset of the set of twisting
cochains $M(B,H_*(G))=\{h=h^2+h^3+...,\ h^k\in
C^k(B,Hom(H_*(G),H_*(G)),\ \delta h=h\smile h\}$. Note that since
the modules $H_i(G)$ are assumed free, an element $h$ is
determined by a collection $\{h(a)\in C^*(B,H_*(G)),\ a\in I\}$
where $I$ is the set of free generators of $H_*(G)$. Now for an
element $(f,h_0)\in L$ we define $h\in M(B,H_*(G))$ as $h=\{h(a),\
a\in I,\ h(a)=h_0\cdot a+\sum_i(\mu f^i)(h_0,...,h_0,a)$.
Inspection shows that $\delta h=h\smile h$. Finally we define
$\alpha(f,h_0)$ as the class of $h$ in $D(B,H_*(G))$.

\vspace{5mm}
 Now let us consider the set of triples
$\bar{L}=\{(f,h_0,\bar{f}) \}$, where $(f,h_0)\in L$ and
$\bar{f}=\bar{f}^3+\bar{f}^4+..., \ \bar{f}^i:H_*(G)\otimes
...((i-1)- times)...\otimes H_*(G)\otimes H_*(F)\to H(F)$ such
that the following condition is satisfied
\begin{equation}
\label{1}
\begin{array}{l}
 \sum_{s+t=i+1}
\bar{f}^s(a_1,...,a_k,f^t(a_{k+1},...,a_{k+t}),...,a_{i-1},x)+\\
\sum_{s+t=i+1}
\bar{f}^s(a_1,...,a_{i-t},\bar{f}^t(a_{i-t+1},...,a_{i-1},x))=\delta
\bar{f}^i(a_,...,a_{i-1},x).
\end{array}
\end{equation}

\vspace{5mm} As above each polylinear map $\bar{f}^i$ induces the
map
$$
(\bar{\mu}\bar{f}^i):C^*(B,H_*(G))\otimes ...\otimes
C^*(B,H_*(G))\otimes C^*(B,H_*(F))\to C^*(B,H_*(F)).
$$

{\bf Definition of the map $\beta:\bar{L}\to D^*(B,H_*(F))$.} We
define $\beta (f,f_0,\bar{f})$ as the class of the twisting
cochain $\bar{h}\in
M(B,H_*(F))=\{\bar{h}=\bar{h}^2+\bar{h}^3+...,\ \bar{h}^k\in
C^k(B,Hom(H_*(F),H_*(F)),\ \delta h=h\smile h\}$, defined as
follows: $\bar{h}=\{\bar{h}(x);\ x\in J\}$ (here $J$ is the set of
free generators of $H_*(F)$), where $\bar{h}(x)=h_0\cdot
x+\sum_i(\bar{\mu}\bar{f}^i)(h_0,...,h_0,x)$. The condition
(\ref{1}) allows to check that $\delta \bar{h}=\bar{h}\smile
\bar{h}$.

\vspace{5mm}
 The group $G$ and, respectively, the action
$G\times F\to F$, define (non uniquely) certain Hochschild
twisting cochain $f$, and respectively the element $\bar{f}$, for
which the condition (\ref{1}) is satisfied. This can be done as
follows.

Since the groups $H_i(G)$ are assumed free, it is possible to fix
a {\it cycle choosing} homomorphism $g:H_*(G)\to Z_*(G)$. Besides
let us fix also a homomorphism $\delta^{-1}:B_*(G)\to C_*(G)$
which satisfies $\delta\delta^{-1}=id$. Let us define also a
homomorphism $\phi:Z_*(G)\to C_*(G)$ by
$\phi(z)=\delta^{-1}(z-g(cl(z)))$, where $cl(z)$ is the homology
class of $z$..

We construct by induction a sequences of multioperations
$f^i=f^3+f^4+...$ and homomorphisms $A_i:H_*(G)\otimes ...\otimes
H_*(G)\to C_*(G)$ using the following conditions:
$$
\begin{array}{l}
1)\ A_2(a_1,a_2)=g(a_1)\cdot g(a_2); \\
2)\ A_i(a_1,...,a_i)\in Z_*(G);\\
3)\ f^i(a_1,...,a_i)=cl(A_i(a_1,...,a_i))\in H_*(G);\\
4)\ A_{i+1}(a_1,...,a_{i+1})=g(a_1)\cdot \phi
A_i(a_2,...,a_{i+1})+\phi A_i(a_1,...,a_{i})\cdot g(a_{i+1})+\\
\sum_{s,t,k} \phi
A_s(a_1,...,a_k,f^t(a_{k+1},...,a_{i+t}),...,a_{i+1})+ \\
 \sum_{s,t}
\phi A_s(a_1,...,a_{s})\cdot \phi A_t(a_{s+1},...,a_{i+1}).
 \end{array}
$$

 Now we construct $\bar{f}$. Let $\bar{g}:H_*(F)\to Z_*(F)$ be a
{\it cycle choosing} homomorphism; $\delta^{-1}:B_*(F)\to C_*(F)$,
$\delta\delta^{-1}=id$; and $\psi:Z_*(F)\to C_*(F)$ is given by
$\psi(z)=\delta^{-1}(z-\bar{g}(cl(z)))$. As above we construct by
induction a sequences of multioperations
$\bar{f}^i=\bar{f}^3+\bar{f}^4+...$ and homomorphisms
$\bar{A}_i:H_*(G)\otimes ...\otimes H_*(G)\otimes H_*(F)\to
C_*(F)$ using the following conditions:
$$
\begin{array}{l}
1)\ \bar{A}_2(a,x)=g(a)\cdot \bar{g}(x); \\
2)\ \bar{A}_i(a_1,...,a_{i-1},x)\in Z_*(F);\\
3)\ \bar{f}^i(a_1,...,a_{i-1},x)=class(\bar{A}^i(a_1,...,a_{i-1},x))\in H_*(F);\\
4)\ \bar{A}_{i+1}(a_1,...,a_{i},x)=g(a_1)\cdot \psi
\bar{A}^i(a_2,...,a_{i-1},x)+\phi A_i(a_1,...,a_{i})\cdot \bar{g}(x)+\\
\sum_{s,t,k} \psi
\bar{A}_s(a_1,...,a_k,f^t(a_{k+1},...,a_{i+t}),...,a_{i},x)+ \\
 \sum_{s,t}
\psi \bar{A}_s(a_1,...,a_{i-t},\bar{f}^t(a_{i-t+1},...,a_i,x))+\\
 \sum_{s,t}
\phi A_s(a_1,...,a_{s})\cdot \psi \bar{A}_t(a_{s+1},...,a_{i},x).
 \end{array}
$$
\noindent {\bf Theorem.} For the above constructed $f$ and
$\bar{f}$ there exists a cochain $h_0=h_0^2+h_0^3+...,\ h_0^i\in
C^i(B,H_{i-1}(G))$ such that $(f,h_0,\bar{f})\in \bar{L}$,
$\alpha(f,h_0)$ is the predifferential of of the principal bundle
$\xi$, and $\beta(f,h_0,\bar{f})$ is the predifferential of of the
associated bundle $\eta$.

In particular if the group $G$ is $(n-1)$-connected, we obtain the
result of Shih: $h_0^r=0$ for $r<n$, and consequently $h^r=0$; as
for $n\leq r\leq 2n-2$ we have $h^r(a)=h_0^r\cdot a$, and $h_0^r$
are cocycles.

\vspace{10mm}

\noindent A. Razmadze Mathematical Institute \newline    \noindent
of the Georgian Academy of Sciences

\noindent kade@rmi.acnet.ge


\begin{thebibliography}{99}
\bibitem{shih} W. Shih, Inst. Hautes et Sci., Publ. Math. Fr., 13,
1962.

\bibitem{fadel} E. Fadel, W. Hurewicz, Ann. of Math., 68,2, 1958,
314-347.

\bibitem{ber1} N. Berikashvili, Bull. Georgian Acad. Sci., 51, 1, 1968.

\bibitem{ber2} N. Berikashvili, On differentials of spectrel
sequences, Tbilisi, 1971.

\end{thebibliography}
\end{document}